\documentclass[a4paper,10pt]{amsart}
\usepackage[english]{babel}
\usepackage{amssymb,amsmath,amscd,a4,amsfonts,amsthm,fullpage}
\date{\today}
\newtheoremstyle{style1}{10pt}{10pt}{\itshape}{}{\bfseries}{.}%
{7pt}{\thmnumber{\textbf{#2.}\;\;}\thmname{#1}\thmnote{#3}}
\newtheoremstyle{style2}{10pt}{10pt}{}{}{\bfseries}{.}%
{7pt}{\thmnumber{\textbf{#2.}\;\;}\thmname{#1}\thmnote{#3}}
\theoremstyle{style2}\newtheorem{dfn}{Definition}[section]
\theoremstyle{style2}
\theoremstyle{style2}\newtheorem{exms}[dfn]{Examples}
\theoremstyle{style2}\newtheorem{rem}[dfn]{Remark}
\theoremstyle{style2}\newtheorem{rems}[dfn]{Remarks}
\theoremstyle{style2}\newtheorem{nots}[dfn]{Notations}
\theoremstyle{style1}\newtheorem{lem}[dfn]{Lemma}
\theoremstyle{style1}\newtheorem{prop}[dfn]{Proposition}
\theoremstyle{style1}\newtheorem{thm}[dfn]{Theorem}
\theoremstyle{style1}\newtheorem{cor}[dfn]{Corollary}

\begin{document}

\title{On Ozawa kernels}

\author{Ghislain Jaudon}
\address{Universit\'{e} de Gen\`{e}ve,
Section de Math\'{e}matiques, 2-4 rue du Li\`{e}vre, Case postale
64, 1211 Gen\`{e}ve 4, Switzerland}
\email{ghislain.jaudon@math.unige.ch}

\begin{abstract}
\noindent We write explicitly Ozawa kernels for group extensions,
for discrete metric spaces of finite asymptotic
dimension, of large enough Hilbert space compression, and for
suitable actions of countable groups on metric spaces. We also
obtain an alternative proof of stability results concerning Yu's
property A.
\end{abstract}

\thanks{This work was supported by the
Swiss National Science Foundation Grant $\sharp$~PP002-68627.}

\maketitle

\setcounter{section}{0}

\section{Introduction}

\noindent Property A is a weak form of amenability which was first
introduced by G. Yu in \cite{yu}. He builds on ideas developed by
M.B. Bekka, P-A. Cherix and A. Valette in \cite{bcv} and shows that
a finitely generated group with property A is uniformly embeddable
in some Hilbert space, the main result of \cite{yu} being that such
a group satisfies the coarse Baum-Connes
Conjecture and the Novikov Higher Signature Conjecture.\\

\noindent The class of finitely generated groups satisfying property
A contains, for instance, amenable groups, hyperbolic groups (and
groups which are hyperbolic relatively to a finite family of
subgroups with property A), one relator groups, Coxeter groups
(moreover any group of finite asymptotic dimension), any discrete
subgroup of a connected Lie group, groups acting properly on finite
dimensional CAT(0) cube complexes, or more generally every group
acting by isometries on a metric space (with bounded geometry)
having property A with at least one point stabilizer having property
A. Furthermore, Yu's property A is known to be closed under taking
subgroups, extensions, direct limits, amalgamated
free products and HNN extensions (see for instance \cite{tu},
\cite{cn}, \cite{dg}, \cite{gu}, \cite{ghw} and \cite{hr}). Actually, the only
known examples of groups which do not satisfy
property A are due to M. Gromov (see \cite{gr}).\\

\noindent Property A admits several equivalent definitions. Here we
focus on a formulation of this property in terms of the existence of
an approximation of the unity by positive definite kernels of finite
width (called Ozawa kernels, see definitions below). One aim of this
paper is to study the behaviors of these kernels and to find explicit
formulas. We write explicitly these kernels for group extensions, for discrete metric spaces
of finite asymptotic dimension, for discrete metric spaces of large
enough Hilbert space compression and moreover for groups acting in a
suitable way on metric spaces with property A. In the last section we
apply formulas obtained to particular examples like hyperbolic groups, CAT(0) cubical groups
and Baumslag-Solitar groups.\\

\noindent\textbf{Acknowledgements.}\\
\noindent The author would like to thank Goulnara Arzhantseva,
Pierre-Alain Cherix, Erik Guentner, Pierre De La Harpe, Graham Niblo, Narutaka Ozawa,
Hugo Parlier and Alain Valette for helpful discussions and valuable remarks. The author also thanks the reviewer for his comments.

\section{Property A and Ozawa kernels}

\noindent First of all we recall the definition of Yu's property A:

\begin{dfn}\label{definition de A}
A discrete metric space $(X,d)$ is said to have \textit{property A}
if for every $R>0$ and every $\varepsilon>0$ there exists a family
of finite sets $\{A_x\}_{x\in X}$ in $X\times\mathbb{N}$ satisfying:
\begin{enumerate}

\item[(1)] $\exists S>0$ such that $d(x,y)\leq S$ whenever $(x,m)\in A_y$;

\item[(2)] $\forall x,y\in X$ such that $d(x,y)\leq R$, we have $|A_x\triangle A_y| < \varepsilon |A_x\cap A_y|$
($|A|$ denoting the cardinality of $A$).\\
\end{enumerate}
\end{dfn}

\noindent If $X=\Gamma$ is a countable group, up to coarse
equivalence, there is a unique way to endow $\Gamma$ with a left
invariant metric (induced naturally by a proper length function $l$,
i.e., $d(x,y):=l(x^{-1}y)$) for which the resulting metric space has
bounded geometry (see \cite{tu} Lemma 2.1 and Lemma 4.1). In the sequel,
all groups will be considered endowed with such a metric. For finitely generated
groups, the metric will always considered to be induced by a length
function associated to a fixed finite generating set. This metric
will be denoted by ``$d_{\Gamma}$", and $B_{\Gamma}(\gamma,S)$ will
denote the closed ball centered at $\gamma$ of radius
$S$ in $\Gamma$ with respect to that metric.

\begin{thm}\label{equivalence des definitions de A}
Let $\Gamma$ be a countable group, then the following assertions are equivalent:
\begin{enumerate}
\item[$(i)$] $\Gamma$ has property A;

\item[$(ii)$] The action of $\Gamma$ on its Stone-$\check{\textrm{C}}$ech
compactification by left translations is topologically amenable;

\item[$(iii)$] The Roe $C^{\ast}$-algebra of $\Gamma$ is nuclear;

\item[$(iv)$] The reduced $C^{\ast}$-algebra of $\Gamma$ is exact.
\end{enumerate}
\end{thm}

\noindent The equivalence ``$(i)\Leftrightarrow (ii)$" is due to N.
Higson and J. Roe \cite{hr}, the equivalence ``$(ii)\Leftrightarrow
(iii)$" can be found in \cite{delaroche} and the equivalence
``$(iii)\Leftrightarrow (iv)$" is due to N. Ozawa \cite{oz}. In his
proof, N. Ozawa introduces positive definite kernels to emphasize
the links between the geometric properties of a group and properties
of its reduced $C^{\ast}$-algebra.

\begin{dfn}\label{Ozawa}
Let $X$ be a set. A function $\psi:X\times X\rightarrow\mathbb{R}$
is said to be a {\it positive definite kernel} if
$\psi(x,y)=\psi(y,x)$ for all $x,y\in X$, and if for every integer
$n\geq 1$, for every $x_1,\ldots,x_n\in X$ and for every
$\lambda_1,\ldots,\lambda_n\in\mathbb{R}$, the following inequality
holds:
$$
\sum_{1\leq i,j\leq n}\lambda_i\lambda_j\psi(x_i,x_j)\geq 0.
$$
\end{dfn}

\begin{dfn}
A discrete metric space $X$ satisfies \textit{Ozawa's property} if
for every $R>0$ and every $\varepsilon>0$ there exist a positive
definite kernel $\psi:X\times X\rightarrow\mathbb{R}$ and a constant
$S\geq R$ such that $\textrm{supp}(\psi)\subset\{(x,y)\in X\times X
\mid d(x,y)\leq S\}$ and $|1 - \psi(x,y)|<\varepsilon$ for every
$x,y\in X$ such that $d(x,y)\leq R$. Such kernels $\psi$ are called
\textit{Ozawa kernels} (or, more precisely,
\textit{$(R,\epsilon)$-Ozawa kernels}).
\end{dfn}

\noindent In the case of countable groups, the following result is a
consequence of Theorem $\ref{equivalence des definitions de A}$, but
a direct proof (without any reference to $C^{\ast}$-algebras) is the
subject of Proposition 3.2 in \cite{tu} together with Lemma 3.5 of
\cite{hr}:

\begin{thm}\label{A equiv O}
Let $X$ be a discrete metric space of bounded geometry. The
following assertions are equivalent:
\begin{enumerate}
\item[$(i)$] X has property A;

\item[$(ii)$] X has Ozawa's property.
\end{enumerate}
\end{thm}

\begin{rems}
If one wants to compare property A and amenability, Ozawa's
property plays the role of Hulanicki's property (see for instance
appendix G in \cite{bhv}). The proof of J-L. Tu in \cite{tu} allows
us to systematically deduce Ozawa kernels from the sets in
Definition $\ref{definition de A}$.
\end{rems}
\begin{prop}\label{noyaux d'ozawa}
Let $X$ be a discrete metric space. Let $\varepsilon>0$, $R>0$, $S>0$, and $\{A_x\}_{x\in X}$ as in Definition $\ref{definition de A}$. Then
$$
\psi_X(x,y):=\sum_{z\in X}\left(\frac{\sum_{n\in\mathbb{N}}\text{{\large$\chi$}}_{A_x}(z,n)}{|A_x|}\right)^{1/2}
\left(\frac{\sum_{n\in\mathbb{N}}\text{{\large$\chi$}}_{A_y}(z,n)}{|A_y|}\right)^{1/2}
$$
is an $(R,2\varepsilon)$-Ozawa kernel, and its support is contained
in $\{(x,y)\in X\times X \mid d(x,y)\leq 2S\}$.
\end{prop}

\begin{proof}[Proof of Proposition $\ref{noyaux d'ozawa}$]
One has $\frac{|A_x\triangle A_y|}{|A_x|}<\varepsilon$, whenever
$d(x,y)\leq R$. This is equivalent to
$\frac{\|\text{{\normalsize$\chi$}}_{A_x} -
\text{{\normalsize$\chi$}}_{A_y}\|}{\|\text{{\normalsize$\chi$}}_{A_x}\|}<
\varepsilon^{1/2}$, whenever $d(x,y)\leq R$ (where
$\|\cdot\|:=\|\cdot\|_{l^2(X\times \mathbb{N})}$). Then,
$$
 \left\|\frac{\text{{\large$\chi$}}_{A_x}}{\|\text{{\large$\chi$}}_{A_x}\|}
 - \frac{\text{{\large$\chi$}}_{A_y}}{\|\text{{\large$\chi$}}_{A_y}\|}\right\|\leq2
 \frac{\|\text{{\large$\chi$}}_{A_x} -
\text{{\large$\chi$}}_{A_y}\|}{\|\text{{\large$\chi$}}_{A_x}\|}<2\varepsilon^{1/2},
$$
whenever $d(x,y)\leq R$. For every $x\in X$, let us consider the
function $\eta_x:X\rightarrow\mathbb{R},z\mapsto
\frac{\|\text{{\normalsize$\chi$}}_{A_x}(z,\cdot)\|_{l^2(\mathbb{N})}}
{\|\text{{\normalsize$\chi$}}_{A_x}\|}=\left(\frac{\sum_{n\in\mathbb{N}}
\text{{\normalsize$\chi$}}_{A_x}(z,n)}{|A_x|}\right)^{1/2}$. We have
$\|\eta_x-\eta_y\|_{l^2(X)}\leq
 \left\|\frac{\text{{\normalsize$\chi$}}_{A_x}}{\|\text{{\normalsize$\chi$}}_{A_x}\|}
 - \frac{\text{{\normalsize$\chi$}}_{A_y}}{\|\text{{\normalsize$\chi$}}_{A_y}\|}\right\|
 < 2\varepsilon^{1/2}$ whenever $d(x,y)\leq R$.\\
Finally, put
$$
\psi_X(x,y):=\langle\eta_x,\eta_y\rangle=\sum_{z\in X}\left(\frac{\sum_{n\in\mathbb{N}}
\text{{\large$\chi$}}_{A_x}(z,n)}{|A_x|}\right)^{1/2}
\left(\frac{\sum_{n\in\mathbb{N}}\text{{\large$\chi$}}_{A_y}(z,n)}{|A_y|}\right)^{1/2},
$$
This is a positive definite kernel satisfying
$\textrm{supp}(\psi_X)\subset\{(x,y)\in X\times X \mid d(x,y)\leq
2S\}$ and $1- \psi_X(x,y)=\frac{1}{2}\|\eta_x -\eta_y\|_{l^2(X)}^2 <
2\varepsilon$ for all $x,y\in X$ with $d(x,y)\leq R$.
\end{proof}

\begin{rem}

If there is only one ``level" in the definition of property A, i.e.,
if $X\times\mathbb{N}$ may be replaced by $X$ in Definition
$\ref{definition de A}$, then Ozawa kernels can be taken of the
form $(x,y)\mapsto \frac{|A_x\cap A_y|}{\sqrt{|A_x||A_y|}}=
\frac{1}{\sqrt{|A_x||A_y|}}\sum_{z\in
X}\text{{\large$\chi$}}_{A_x}(z)\text{{\large$\chi$}}_{A_y}(z)$.
This is actually the case for all known examples of groups with
property A. In particular, using examples 2.3 and 2.4 in
\cite{yu}, we easily recover Ozawa kernels obtained by S. Campbell
in \cite{c} (see the following examples).
\\\\
Suppose $X=\Gamma$ is a group, and the Ozawa kernels are
$\Gamma$-invariant. Set $\phi(\gamma):=\psi(e,\gamma)$. Then
we obtain an approximation of the unity in $\Gamma$ by positive
definite \textit{functions} with finite support. Actually this
characterizes amenability: if a group $\Gamma$ has property
A and admits $\Gamma$-invariant Ozawa kernels, then $\Gamma$ is
amenable. The converse is also true, see the first example below.

\end{rem}

\begin{exms}\label{exemples de noyaux}
Let $\varepsilon, R>0$.
\begin{enumerate}

\item Let $\Gamma$ be a countable amenable group. The length function on $\Gamma$ being assumed to be proper, $F:=B_{\Gamma}(e,R)$ is a finite set.
Using Folner's definition of amenability, there exists a finite
subset $A\subset \Gamma$ such that $\frac{|\gamma A\triangle
A|}{|A|} < \frac{2\epsilon}{4 + \epsilon}$ for every $\gamma\in F$.
Then we define $A_{\gamma}:=(\gamma A)\times\{1\}$ (note that
$|A_{\gamma}|=|A|$ for every $\gamma$). If $\gamma,\gamma'\in\Gamma$
are such that $d_{\Gamma}(\gamma,\gamma')\leq R$, we have
$\gamma^{-1}\gamma'\in F$ and $|A_{\gamma}\triangle A_{\gamma'}|=
|(\gamma^{-1}\gamma' A) \triangle A| < \frac{2\epsilon}{4 +
\epsilon}|A_{\gamma}|$. We also have
$$
|A_{\gamma}\cap A_{\gamma'}|=\frac{1}{2}[|A_{\gamma}| + |A_{\gamma'}| - |A_{\gamma}\triangle A_{\gamma'}|]>  \frac{4}{4 + \epsilon}|A_{\gamma}|
$$
Therefore $|A_{\gamma}\triangle A_{\gamma'}|< \frac{\epsilon}{2}|A_{\gamma}\cap A_{\gamma'}|$.
Moreover, if $S:=\max_{a\in A}d_{\Gamma}(e,a)$, $(\gamma,1)\in A_{\gamma'}$, implies $d_{\Gamma}(\gamma,\gamma')\leq S$.

\noindent Hence, the family $\{A_{\gamma}\}_{\gamma\in \Gamma}$ satisfies
Definition $\ref{definition de A}$, and by Proposition $\ref{noyaux
d'ozawa}$, we have shown that
$$
\psi:(\gamma,\gamma')\mapsto \frac{|(\gamma A)\cap (\gamma' A)|}{|A|}
$$
is an $(R,\epsilon)$-Ozawa kernel for $\Gamma$ (and this kernel is
clearly $\Gamma$-invariant).

\item Let $\Gamma$ be a finitely generated free group, and $T$ be the standard Cayley graph (which is a tree) of
$\Gamma$. We fix a geodesic ray $r_0$ in $T$, and for any $\gamma\in
\Gamma$, we denote by $r(\gamma)$ the unique geodesic ray starting
from $\gamma$ such that $r(\gamma)\cap r_0$ is a non-empty geodesic
ray. Let $S > \frac{R}{2}$ such that $\frac{2R}{2(S+1) - R} <
\frac{\epsilon}{2}$, and for every $\gamma\in \Gamma$, let us define
$A_{\gamma}:=[r(\gamma)\cap B_{\Gamma}(\gamma, S)]\times\{1\}\subset
\Gamma\times\mathbb{N}$ (note that $|A_{\gamma}|=S+1$ for every
$\gamma$). Now, let $\gamma,\gamma'\in\Gamma$ such that
$d_{\Gamma}(\gamma,\gamma')\leq R$. As $S>\frac{R}{2}$,
$A_{\gamma}\cap A_{\gamma'}\neq \emptyset$, and
$|A_{\gamma}\triangle A_{\gamma'}|$ is maximal (being equal to
$d_{\Gamma}(\gamma,\gamma')$) exactly when $A_{\gamma}$ and
$A_{\gamma'}$ intersect in just one point, (i.e. when $A_{\gamma}\cup
A_{\gamma'}$ realize the unique geodesic path between $\gamma$ and
$\gamma'$ in $T$), then for any $\gamma,\gamma'\in\Gamma$ such that
$d_{\Gamma}(\gamma,\gamma')\leq R$, one has $|A_{\gamma}\triangle
A_{\gamma'}|\leq R$. We also have
$$
|A_{\gamma}\cap A_{\gamma'}|=\frac{1}{2}[|A_{\gamma}| +
|A_{\gamma'}| - |A_{\gamma}\triangle A_{\gamma'}|]\geq  \frac{2(S+1)
- R}{2}.
$$
Therefore, by our choice of $S$, we obtain
$$
|A_{\gamma}\triangle A_{\gamma'}| \leq  \frac{2R}{2(S+1) -
R}|A_{\gamma}\cap A_{\gamma'}| < \frac{\epsilon}{2}|A_{\gamma}\cap
A_{\gamma'}|.
$$
Moreover, by construction, $(\gamma,1)\in A_{\gamma'}$ implies
$d_{\Gamma}(\gamma,\gamma')\leq S$. Hence, we have shown that
$$
(\gamma,\gamma')\mapsto \frac{|A_{\gamma}\cap A_{\gamma'}|}{S+1}
$$
is an $(R,\epsilon)$-Ozawa kernel for $\Gamma$.
\end{enumerate}
\end{exms}

\begin{rem}
Using, the so-called GNS construction (see for instance
Proposition 3 of Chapter 5 in \cite{hv}), one knows that any Ozawa
kernel on a discrete metric space $X$ is of the form
$\psi:(x,y)\mapsto \langle\lambda(x),\lambda(y)\rangle$, for some
$\lambda:X\rightarrow l^2(X)_1:=\{\xi\in
l^2(X)~\mid~\|\xi\|_{l^2(X)}=1\}$. Unfortunately, in general, we
have no information about the support of the functions in the range
of $\lambda$. However, if $X$ has bounded geometry, as a corollary
of Theorem $\ref{A equiv O}$, we have the following:
\end{rem}

\begin{prop}\label{forme des noyaux}
Let $X$ be a  discrete metric space with bounded geometry, if $X$
has property A, then an $(R,\epsilon)$-Ozawa kernel for $X$
can always be taken of the form $\psi:(x,y)\mapsto
\langle\lambda(x),\lambda(y)\rangle$, for some $\lambda:X\rightarrow
l^2(X)_1$ with $\textrm{supp}(\lambda(x))\subset B(x,S)$ for every
$x$ (and $\lambda(x)(y)\geq 0$ for every $y$) and some uniform
constant $S$ depending only on $R$ and $\epsilon$.
\end{prop}

\section{Group extensions}
\noindent The exactness of $C^{\ast}$-algebras is stable under
taking extensions. This is a result due to E. Kirchberg and S.
Wassermann (see \cite{kw}), but
the proof is quite technical and done in the context of
$C^{\ast}$-algebras. G. Yu \cite{yu} indicates how to prove the
stability of property A under taking semi-direct products leaving
details to the reader. Here we give a detailed simple proof of the
stability under general extensions by expliciting Ozawa kernels. It
is inspired by \cite{dg} (see also \cite{delarocherenault}).

\begin{thm}\label{thm extensions}
Let $1\rightarrow H \xrightarrow{i} \Gamma \xrightarrow{\pi} G \rightarrow 1$ be
a short exact sequence of countable groups. If $H$ and $G$ have
property A then $\Gamma$ has property A.
\end{thm}

\noindent As an easy consequence, we obtain the following:

\begin{cor}\label{corollaire thm extensions}
Let $\Gamma$ be a countable group, and let $H$ be a finite index subgroup. Then property A for $H$ implies property A for $\Gamma$.
\end{cor}

\begin{proof}[Proof of Corollary $\ref{corollaire thm extensions}$]
Let $\{\gamma_1,\ldots,\gamma_n\}$ be a set of representatives for the left cosets of $H$ in $\Gamma$. Then $N:=\bigcap_{i=1}^n \gamma_i H \gamma_i^{-1}$ is also a finite index subgroup of $\Gamma$ ($[\Gamma:N]\leq [\Gamma:H]^n$), and moreover (by construction) $N$ is normal in $\Gamma$. Being a subgroup of $H$, $N$ has also property A. Now we can conclude applying theorem $\ref{thm extensions}$ to the extension $1\rightarrow N \rightarrow \Gamma \rightarrow \Gamma/N \rightarrow 1$ using the fact that a finite group obviously has property A.
\end{proof}

\begin{rem}
Actually the preceding result is also true if we replace the finite index hypothesis by the property A of $G/H$ viewed as an abstract metric space endowed with the quotient metric (see Corollary $\ref{groupe ayant sous groupe et quotient avec A}$).
\end{rem}

\begin{nots}\label{notations extensions}
For the proof of theorem \ref{thm extensions}, given a group extension $1\rightarrow H \xrightarrow{i} \Gamma \xrightarrow{\pi} G \rightarrow 1$, we deduce Ozawa kernels
for $\Gamma$ from Ozawa kernels for $H$ and $G$. Let $l_{\Gamma}$ be an arbitrary proper length function on $\Gamma$.
In the sequel, we identify $H$ with $i(H)\subset\Gamma$, and then we consider
the restriction $l_H:=(l_{\Gamma})_{\mid_H}$ of $l_{\Gamma}$ on $H$. On $G$, we will consider the
length function $l_G$ defined by
$$
l_G(g):=\min\{l_{\Gamma}(\gamma) \mid \gamma\in\Gamma,~\pi(\gamma)=g\}.
$$
One can easily check that $l_G$ is a proper length function
and that the minimum is realized. We endow $G$, $H$ and $\Gamma$ with the left invariant
metrics denoted $d_G$, $d_H$ and $d_{\Gamma}$ associated respectively to $l_G$, $l_H$ and $l_{\Gamma}$. Note that $l_G$ and $l_H$ being proper, the associated metrics on $G$ and $H$ are coarsely equivalent to the original metrics on this groups.\\

\noindent Let us fix a set theoretical section $\sigma$ of $\pi$ in such a way
that $l_{\Gamma}(\sigma(g))=l_G(g)$ for every $g\in G$. The idea of the proof is then to decompose
elements of $\Gamma$, using $\sigma$, into two parts, one in $H$ and one in $G$ as in the case
of a split sequence. Hence we will need the following cocycle:
$$
c:\Gamma\times G\rightarrow H,~(\gamma,g)\mapsto \sigma(g)^{-1}\gamma\sigma(\pi(\gamma)^{-1}g).
$$
\end{nots}
\noindent For any short exact sequence, the map
$\gamma\mapsto (\pi(\gamma),c(\gamma,\pi(\gamma)))$ is a bijection between $\Gamma$ and $G\times H$.\\

\noindent The following inequalities will be useful.

\begin{lem}\label{lem distances}
For any $\gamma_1,\gamma_2\in\Gamma$, and any $g\in G$, one has:
\begin{enumerate}
\item[$(i)$] $d_G(\pi(\gamma_1),\pi(\gamma_2))\leq d_{\Gamma}(\gamma_1,\gamma_2)$;

\item[$(ii)$] $d_{\Gamma}(\gamma_1, \gamma_2)\leq d_G(g,\pi(\gamma_1))  +
d_H(c(\gamma_1,g),c(\gamma_2,g))+ d_G(g,\pi(\gamma_2))$;

\item[$(iii)$]  $d_H(c(\gamma_1,g),c(\gamma_2,g))\leq d_G(g,\pi(\gamma_1)) + d_{\Gamma}(\gamma_1,\gamma_2)
 + d_G(g,\pi(\gamma_2))$.
\end{enumerate}
\end{lem}

\begin{proof}[Proof of Lemma $\ref{lem distances}$] Inequality $(i)$ comes from the definition of $l_G$. To see inequalities $(ii)$ and $(iii)$,
it suffices to use sub-additivity of length functions and to write that:
$$
\gamma_1^{-1}\gamma_2=\sigma(\pi(\gamma_1)^{-1}g)c(\gamma_1,g)^{-1}c(\gamma_2,g)(\sigma(\pi(\gamma_2)^{-1}g))^{-1}.
$$

\end{proof}

\noindent By Theorem $\ref{A equiv O}$, Theorem $\ref{thm extensions}$ is a consequence of the following statement:

\begin{thm}\label{noyaux d'ozawa pour extensions}
Let $1\rightarrow H \xrightarrow{i} \Gamma \xrightarrow{\pi} G \rightarrow 1$ be
a short exact sequence of countable groups, $H$ and $G$ having property A. Let
$\psi_G:(g_1,g_2)\mapsto \langle\lambda(g_1),\lambda(g_2)\rangle$ and  $\psi_H:(h_1,h_2)\mapsto \langle\mu(h_1),\mu(h_2)\rangle$
be Ozawa kernels for $G$ and $H$ given by Proposition $\ref{forme des noyaux}$ (with $\lambda:G\rightarrow l^2(G)_1$ and $\mu:H\rightarrow l^2(H)_1$).
Then maps of the form
$$
\psi_{\Gamma}:(\gamma_1,\gamma_2)\mapsto\sum_{g\in G}\lambda(\pi(\gamma_1))(g)\cdot\lambda(\pi(\gamma_2))(g)\cdot\psi_H(c(\gamma_1,g),c(\gamma_2,g))
$$
are Ozawa kernels for $\Gamma$ (where $c$ denotes the cocycle defined above).
\end{thm}

\begin{proof}[Proof of Theorem $\ref{noyaux d'ozawa pour extensions}$]
First, note that, because $\lambda(\pi(\gamma))$ has finite support, the sum defining $\psi_{\Gamma}(\gamma_1,\gamma_2)$ is a finite sum.
To check that $\psi_{\Gamma}$ is indeed a positive definite kernel, it suffices to note that
for every $(\gamma_1,\gamma_2)\in \Gamma\times\Gamma$,
$\psi_{\Gamma}(\gamma_1,\gamma_2)=\langle\nu(\gamma_1),\nu(\gamma_2)\rangle$, where
$$
\nu:\Gamma\rightarrow l^2(G,l^2(H))~,~\gamma\mapsto\nu(\gamma):g\mapsto
\lambda(\pi(\gamma))(g)\cdot\mu(c(\gamma,g)).
$$

\noindent Now, let $\varepsilon>0$ and $R>0$ be fixed. By hypothesis,
there exist positive constants $S_1$ and $S_2$
such that:\\

\begin{enumerate}
\item[$(a)$] $d_G(g_1,g_2)\leq R \Rightarrow |1 - \psi_G(g_1,g_2)| < \frac{\varepsilon}{2}$;
\item[$(b)$] $d_H(h_1,h_2)\leq R + 2S_1 \Rightarrow |1 - \psi_H(h_1,h_2)| < \frac{\varepsilon}{2}$;
\item[$(c)$] $\textrm{supp}(\lambda(g))\subseteq B_G(g,S_1)$, $\forall g\in G$;
\item[$(d)$] $d_H(h_1,h_2)> S_2\Rightarrow \psi_H(h_1,h_2)=0$.
\end{enumerate}

\noindent For any $\gamma_1,\gamma_2\in\Gamma$, we have:
$$
|1 - \psi_{\Gamma}(\gamma_1,\gamma_2)|\leq \left|\sum_{g\in G}\left(1 - \psi_H(c(\gamma_1,g),c(\gamma_2,g))\right)\cdot
\lambda(\pi(\gamma_1))(g)\cdot\lambda(\pi(\gamma_2))(g)\right|
 + ~|1 - \psi_G(\pi(\gamma_1),\pi(\gamma_2))|.~~~~~(\star)
$$
Then, if $d_{\Gamma}(\gamma_1,\gamma_2)\leq R$, by Lemma $\ref{lem distances}$ $(i)$ and by $(a)$, one has
$|1 - \psi_G(\pi(\gamma_1),\pi(\gamma_2))| < \frac{\varepsilon}{2}$.\\

\noindent On the other hand, by $(c)$, the sum in the first term of $(\star)$
is just over the set of $g\in G$ such that $d_G(g,\pi(\gamma_1))\leq S_1$ and $d_G(g,\pi(\gamma_2))\leq S_1$. But for such $g\in G$,
by $(iii)$ of Lemma $\ref{lem distances}$, one has $d_H(c(\gamma_1,g),c(\gamma_2,g))\leq R + 2S_1$. Therefore, by $(b)$ and by the Cauchy-Schwarz inequality\\
(as $\sum_{g\in G}\lambda(x)(g)^2=1$ for all $x\in G$), one obtains:
$$
\left|\sum_{g\in G}\left(1 - \psi_H(c(\gamma_1,g),c(\gamma_2,g))\right)\cdot
\lambda(\pi(\gamma_1))(g)\cdot\lambda(\pi(\gamma_2))(g)\right|< \frac{\varepsilon}{2}.
$$
Hence, $|1 - \psi_{\Gamma}(\gamma_1,\gamma_2)| < \varepsilon$ for all $\gamma_1,\gamma_2\in\Gamma$ such that $d_{\Gamma}(\gamma_1,\gamma_2)\leq R$.\\

\noindent Moreover, again by $(c)$ and by the Cauchy-Schwarz inequality, we have:
$$
|\psi_{\Gamma}(\gamma_1,\gamma_2)|\leq\max\{|\psi_H(c(\gamma_1,g),c(\gamma_2,g))| ~\mid g\in G~,~d_G(g,\pi(\gamma_k)\leq S_1~,~k=1,2\}.~~~~~~~~~~~~~~~~~~~~~~~~~~~~~~~~~(\star\star)
$$

\noindent If $\gamma_1,\gamma_2\in\Gamma$ are that $d_{\Gamma}(\gamma_1,\gamma_2)> 2S_1 + S_2$, by $(ii)$ of Lemma $\ref{lem distances}$, one has:
$$
d_H(c(\gamma_1,g),c(\gamma_2,g))\geq d_{\Gamma}(\gamma_1,\gamma_2) - d_G(g,\pi(\gamma_1)) - d_G(g,\pi(\gamma_2))> S_2.
$$
Thus, by $(d)$, we deduce from $(\star\star)$ that
$\textrm{supp}(\psi_{\Gamma})\subset\{(\gamma_1,\gamma_2)\in \Gamma\times\Gamma \mid d_{\Gamma}(\gamma_1,\gamma_2)\leq 2S_1 + S_2\}$.
This concludes the proof.

\end{proof}

\begin{rem}\label{remarque sur les cocycles}
Suppose that short exact sequence splits, i.e., when $\Gamma\simeq H\rtimes G$. Then viewing $G$ as a subgroup of $\Gamma$, the metric corresponding to the restriction to $G$ of $l_{\Gamma}$ is coarsely equivalent to the metric $d_G$ defined above. Hence, in this case, up to coarse equivalence, $c(\gamma,g)$ may be replaced by
the action of $g^{-1}\pi(\gamma)$ on  $\pi(\gamma)^{-1}\gamma$ (the component of $\gamma$ in $H$).
\end{rem}

\section{Metric spaces with finite asymptotic dimension}

\noindent Among finitely generated groups, there are important examples of groups of finite asymptotic dimension, for instance, Coxeter groups \cite{dj}, one relator groups \cite{m} and hyperbolic groups \cite{roe}. It is well known that for a discrete metric space $X$ of bounded geometry, finite asymptotic dimension implies property A (see \cite{hr} and \cite{dg1}). We follow the proof of this result but make it more precise as far as Ozawa kernels are concerned. First, recall:

\begin{dfn}\label{asdim}
Let $(X,d)$ be a metric space and $k\in\mathbb{N}$. The space $X$ is said to have \textit{asymptotic dimension less than $k$}, denoted by $\textrm{asdim}X\leq k$, if for any $L>0$, there is an open cover $\mathcal{U}:=\{U_i\}_{i\in I}$ of $X$ satisfying:
\begin{enumerate}
\item[$(i)$] $\sup_{i\in I} \textrm{diam}(U_i) < \infty$ ($\mathcal{U}$ is said to be \textit{uniformly bounded});
\item[$(ii)$] The \textit{Lebesgue number} of $\mathcal{U}$, $L(\mathcal{U}):=\inf\{\max\{d(x,X\smallsetminus U_i)~\mid~i\in I\}~\mid~x\in X\}$, satisfies $L(\mathcal{U})\geq L$;
\item[$(iii)$] For every $x\in X$, $|\{i\in I~\mid~x\in U_i\}|\leq k+1$  ($\mathcal{U}$ is said to have \textit{multiplicity $\leq k+1$}).
\end{enumerate}

\end{dfn}

\begin{lem}\label{noyau de recouvrement}
Let $(X,d)$ be a discrete metric space admitting a uniformly bounded cover $\mathcal{U}:=\{U_i\}_{i\in I}$ with multiplicity $\leq k+1$. Then the kernel $\psi_{\mathcal{U}}$ defined by
$$
\psi_{\mathcal{U}}:X\times X\rightarrow\mathbb{R}~,~(x,y)\mapsto\sum_{i\in I}\left(\frac{d(x,X\smallsetminus U_i)}{\sum_{j\in I}d(x,X\smallsetminus U_j)}\right)^{1/2}\left(\frac{d(y,X\smallsetminus U_i)}{\sum_{j\in I}d(y,X\smallsetminus U_j)}\right)^{1/2}
$$
is positive definite and satisfies $|1-\psi_{\mathcal{U}}(x,y)|\leq \frac{(k+1)(2k+3)}{L(\mathcal{U})}d(x,y)$ for every $x,y\in X$.
\end{lem}

\begin{proof}[Proof of Lemma $\ref{noyau de recouvrement}$]
The kernel $\psi_{\mathcal{U}}$ is a positive definite kernel as $\psi_{\mathcal{U}}(x,y)=\langle\lambda(x),\lambda(y)\rangle$, where $\lambda:X\rightarrow l^2(I)~,~x\mapsto(\frac{d(x,X\smallsetminus U_i)}{\sum_{j\in I}d(x,X\smallsetminus U_j)})_{i\in I}$. Moreover, we have
\begin{align*}
1-\psi_{\mathcal{U}}(x,y) &=\frac{1}{2}\|\lambda(x)-\lambda(y)\|_{l^2(I)}^2=\frac{1}{2}\sum_{i\in I}\left|\left(\frac{d(x,X\smallsetminus U_i)}{\sum_{j\in I}d(x,X\smallsetminus U_j)}\right)^{1/2} - \left(\frac{d(y,X\smallsetminus U_i)}{\sum_{j\in I}d(y,X\smallsetminus U_j)}\right)^{1/2}\right|^2\\
&\leq \frac{1}{2}\sum_{i\in I}\left|\frac{d(x,X\smallsetminus U_i)}{\sum_{j\in I}d(x,X\smallsetminus U_j)} - \frac{d(y,X\smallsetminus U_i)}{\sum_{j\in I}d(y,X\smallsetminus U_j)}\right|\\
&\leq (k+1)\max_{i\in I}\left|\frac{d(x,X\smallsetminus U_i)}{\sum_{j\in I}d(x,X\smallsetminus U_j)} - \frac{d(y,X\smallsetminus U_i)}{\sum_{j\in I}d(y,X\smallsetminus U_j)}\right|
\end{align*}
the last inequality coming from the fact that the multiplicity of $\mathcal{U}$ is less than $k+1$, i.e.,
$$
|\{i\in I~\mid~d(x,X\smallsetminus U_i)~~\textrm{or}~~ d(y,X\smallsetminus U_i)\neq 0\}|\leq 2k+2~~~~~~~~~~~~~~~~~~~~~~~~~~~~~~~~~~~~~~~~~~~~~~~~~~~~~~~~~~~~~~(*)
$$
However for each $i\in I$
\begin{align*}
&\left|\frac{d(x,X\smallsetminus U_i)}{\sum_{j\in I}d(x,X\smallsetminus U_j)} - \frac{d(y,X\smallsetminus U_i)}{\sum_{j\in I}d(y,X\smallsetminus U_j)}\right|\\
&\leq \left|\frac{d(x,X\smallsetminus U_i)}{\sum_{j\in I}d(x,X\smallsetminus U_j)} -\frac{d(y,X\smallsetminus U_i)}{\sum_{j\in I}d(x,X\smallsetminus U_j)}\right| + \left|\frac{d(y,X\smallsetminus U_i)}{\sum_{j\in I}d(x,X\smallsetminus U_j)} - \frac{d(y,X\smallsetminus U_i)}{\sum_{j\in I}d(y,X\smallsetminus U_j)}\right|
\end{align*}
Then, as $\sum_{j\in I}d(x,X\smallsetminus U_j)\geq L(\mathcal{U})$ and $\frac{d(y,X\smallsetminus U_i)}{\sum_{j\in I}d(y,X\smallsetminus U_j)}\leq 1$, we obtain
\begin{align*}
&\left|\frac{d(x,X\smallsetminus U_i)}{\sum_{j\in I}d(x,X\smallsetminus U_j)} - \frac{d(y,X\smallsetminus U_i)}{\sum_{j\in I}d(y,X\smallsetminus U_j)}\right|\\
&\leq \frac{1}{L(\mathcal{U})}\left[\left|d(x,X\smallsetminus U_i) - d(y,X\smallsetminus U_i)\right| + \sum_{j\in I}\left|d(x,X\smallsetminus U_j)-d(y,X\smallsetminus U_j)\right|\right]
\end{align*}
therefore, by $(*)$ and by the inequality $\left|d(x,X\smallsetminus U_j) - d(y,X\smallsetminus U_j)\right|\leq d(x,y)$, we deduce
$$
\left|\frac{d(x,X\smallsetminus U_i)}{\sum_{j\in I}d(x,X\smallsetminus U_j)} - \frac{d(y,X\smallsetminus U_i)}{\sum_{j\in I}d(y,X\smallsetminus U_j)}\right|\leq \frac{(2k+3)R}{L(\mathcal{U})}.
$$
\end{proof}

\begin{thm}\label{noyaux pour groupes avec asdim finie}
Let $(X,d)$ be a discrete metric space with $\textrm{asdim}X\leq k$, and let $\epsilon,R>0$. We fix a uniformly bounded cover of $X$ with multiplicity $\leq k+1$, $\mathcal{U}:=\{U_i\}_{i\in I}$, such that $L(\mathcal{U}) > \frac{(k+1)(2k+3)R}{\epsilon}$. Then the kernel
$\psi_{\mathcal{U}}$ defined in Lemma $\ref{noyau de recouvrement}$ is an $(R,\epsilon)$-Ozawa kernel for $X$. In particular, $X$ has property A.
\end{thm}

\begin{proof}[Proof of Theorem $\ref{noyaux pour groupes avec asdim finie}$]
If $S:=\sup_{i\in I} \textrm{diam}(U_i)$, then for every $x,y\in X$ such that $d(x,y) > S$ and for any $i\in I$, we have $d(x,X\smallsetminus U_i)=0$ or $d(y,X\smallsetminus U_i)=0$, thus $\psi_{\mathcal{U}}(x,y)=0$. Moreover, for every $x,y\in X$ such that $d(x,y)\leq R$, by Lemma $\ref{noyau de recouvrement}$ and by our choice for $L(\mathcal{U})$, we have
$$
|1-\psi_{\mathcal{U}}(x,y)|\leq \frac{(k+1)(2k+3)R}{L(\mathcal{U})}\leq \epsilon.
$$
\end{proof}

\section{Metric spaces with Hilbert space compression $> 1/2$}
\noindent It is well known that any discrete metric space with
bounded geometry having property A is uniformly embeddable in a
Hilbert space. The converse is not known. A result of E. Guentner
and J. Kaminker gives a partial converse: a finitely generated
group with Hilbert space compression $> 1/2$ has property A (see
Theorem 3.2 in \cite{gk}). In this section, we give a proof
(strongly inspired by \cite{gk}) of a slightly more general result.
First of all, let us recall the definition of the Hilbert space
compression.

\begin{dfn}
Let $(X,d)$ be a metric space and $\mathcal{H}$ be an Hilbert space. A map $f:X\rightarrow\mathcal{H}$ is said to be a \textit{uniform embedding} if there exist non-decreasing functions $\rho_{\pm}(f):\mathbb{R}_+\rightarrow\mathbb{R}_+$ such that:
\begin{enumerate}
\item $\rho_-(f)(d(x,y))\leq \|f(x) - f(y)\|_{\mathcal{H}}\leq\rho_+(f)(d(x,y))$, for all $x,y\in X$;
\item $\lim_{r\rightarrow +\infty}\rho_{\pm}(f)(r)=+\infty$.
\end{enumerate}
Then the \textit{Hilbert space compression} of the metric space $X$, denoted by $R(X)$, is defined as the sup of all $\beta\geq 0$ for which there exists a uniform embedding into a Hilbert space $f$ with $\rho_+(f)$ affine and $\rho_-(f)(r)=r^{\beta}$ (for $r$ large enough).
\end{dfn}

\begin{dfn}\label{definition quasi-geodesique}
A metric space $(X,d)$ is said to be \textit{quasi-geodesic} if there exist $\delta>0$ and $\lambda\geq 1$ such that for all $x,y\in X$ there exists a sequence $x_0=x, x_1,\ldots, x_n=y$ of elements of $X$ satisfying $d(x_{i-1},x_i)\leq\delta$ for every $i=1,\ldots, n$, and
$$
\sum_{i=1}^nd(x_{i-1},x_i)\leq \lambda d(x,y).
$$
Such a sequence $x_0=x, x_1,\ldots, x_n=y$ will be called a \textit{$(\lambda,\delta)$-chain of length $n$ from $x$ to $y$}.
\end{dfn}

\noindent We will need the following two lemmas.

\begin{lem}\label{metrique à valeurs entières}
Let $(X,d)$ be a quasi-geodesic metric space. There exists an integer-valued metric $\widetilde{d}$ on $X$ such that $(X,\widetilde{d})$ is quasi-geodesic and quasi-isometric to $(X,d)$.
\end{lem}

\begin{proof}[Proof of Lemma \ref{metrique à valeurs entières}]
Let $\delta>0$ and $\lambda\geq 1$ as in Definition \ref{definition quasi-geodesique}. We define a metric $\widetilde{d}$ on $X$ by setting
$$
\widetilde{d}(x,y):=\min\{n\in\mathbb{N}~\mid~\textrm{there exists a $(\lambda,\delta)$-chain of length $n$ from $x$ to $y$}\}.
$$
This integer-valued metric is clearly quasi-geodesic (with $\lambda=\delta=1$). It remains to show that $(X,\widetilde{d})$ is quasi-isometric to $(X,d)$.
Let us fix $x,y\in X$ and set $\widetilde{d}(x,y):=n$. On the one hand,  by defintion of $n$, one can find a $(\lambda,\delta)$-chain of length $n$ from $x$ to $y$. In particular we have
$$
d(x,y)\leq \sum_{i=1}^nd(x_{i-1},x_i)\leq \delta n=\delta\widetilde{d}(x,y).
$$
On the other hand, again by definition of $n$, we must have $d(x_{i-1},x_i) + d(x_i,x_{i+1}) > \delta$ for every $i=1,\ldots, n-1$ (otherwise, we could forget one of the $x_i$'s and this would contradict the minimality of $n$). Hence, in the chain, there is at least $[\frac{n}{2}]\geq\frac{n-1}{2}$ successive distances which are greater than $\frac{\delta}{2}$. Then, we obtain that
$$
\left(\frac{n-1}{2}\right)\frac{\delta}{2}\leq \sum_{i=1}^nd(x_{i-1},x_i)\leq\lambda d(x,y)
$$
and we deduce that
$$
\widetilde{d}(x,y)\leq \frac{4\lambda}{\delta}d(x,y) + 1.
$$
\end{proof}

\begin{lem}\label{croissance}
Let $(X,d)$ be a discrete quasi-geodesic metric space of bounded geometry. Then the growth of $X$ is at most exponential, i.e., there exist constants $B,L>0$ such that $|B(x,R)|\leq BL^R$ for every $R>0$ and for every $x\in X$ (where $B(x,R):=\{y\in X~\mid~d(x,y)\leq R\}$ denotes the closed ball of radius $R$ centered at $x$).
\end{lem}

\begin{proof}[Proof of Lemma \ref{croissance}]
Let $\delta>0$ and $\lambda\geq 1$ as in Definition \ref{definition quasi-geodesique}. By the bounded geometry assumption on $X$, $B_{\delta}:=\max_{x\in X}|B(x,\delta)| < \infty$. Let us fix $R>0$ and $x\in X$. In order to estimate the number of elements in $B(x,R)$ it suffices to estimate the number of minimal $(\lambda,\delta)$-chains starting at $x$ with end-point in $B(x,R)$. On the one hand, if $x_0=x, x_1,\ldots, x_n=y$ is a $(\lambda,\delta)$-chain of length $n=\widetilde{d}(x,y)$ with $y\in B(x,R)$, by the proof of Lemma \ref{metrique à valeurs entières}, we have $n\leq \frac{4\lambda}{\delta}R + 1$. On the other hand, by the definition of $B_{\delta}$, for a fixed $n$, the number of $(\lambda,\delta)$-chains of length $n$ starting at $x$ with end-point in $B(x,R)$ is at most $B_{\delta}^n$. Therefore, we deduce that
$$
|B(x,R)|\leq \sum_{n=0}^{[\frac{4\lambda R}{\delta}]+1}B_{\delta}^n\leq BL^R
$$
with $L:=B_{\delta}^{\frac{4\lambda}{\delta}}$ and $B:=B_{\delta}^2$.
\end{proof}

\begin{thm}\label{noyaux quand compression > 1/2}
Let $(X,d)$ be a discrete quasi-geodesic metric space of bounded geometry such that $R(X)>1/2$, then $X$ has property A.
\end{thm}

\begin{proof}[Proof of Theorem $\ref{noyaux quand compression > 1/2}$]
By Lemma \ref{metrique à valeurs entières}, we can suppose that the metric $d$ on $X$ takes integer values. For every $n\in\mathbb{N}$ and $x\in X$, we will denote $S(x,n):=\{y\in X~\mid~d(x,y)=n\}$.\\

\noindent Let $R>0$ and $0<\epsilon<1$ fixed. By hypothesis, there exist a Hilbert space $\mathcal{H}$, a map $f:X\rightarrow\mathcal{H}$, and positive constants $n_0, \alpha, C, D$ such that
$\|f(x) - f(y)\|_{\mathcal{H}}\leq Cd(x,y)+D$ for every $x,y\in X$ and
$$
\|f(x) - f(y)\|_{\mathcal{H}}\geq n^{\frac{1+\alpha}{2}}
$$
for every $n,x,y\in X$ such that $d(x,y)\geq n\geq n_0$.\\

\noindent Then, for a fixed $k > \frac{CR+D}{-\ln(1-\frac{\epsilon}{2})}$, we consider
$$
\phi:X\times X\rightarrow\mathbb{R}~,~ (x,y)\mapsto \exp\left(-\frac{\|f(x) - f(y)\|_{\mathcal{H}}^2}{k}\right)
$$
By the Schoenberg's theorem (see for instance \cite{hv}), $\phi$ is a positive definite kernel on $X$. Moreover, by our choice of $k$, for every $x,y\in X$ such that $d(x,y)\leq R$ we have
$$
|1 - \phi(x,y)| < \frac{\epsilon}{2}.
$$
Unfortunately, $\phi$ does not have finite width and therefore is not an Ozawa kernel on $X$. Actually, the remainder of the proof will be devoted to approximate $\phi$ in a suitable way using the hypothesis on the Hilbert space compression of $X$. Let us formally consider the kernel operator associated to $\phi$ defined on $\xi\in l^2(X)$ by
$$
\mathcal{U}_{\phi}(\xi)(x):=\sum_{y\in X}\phi(x,y)\xi(y).
$$
The key point of the proof is that $\mathcal{U}_{\phi}$ defines a bounded positive operator from $l^2(X)$ into itself. Indeed, for every $\xi\in l^2(X)$, on the one hand, by Cauchy-Schwarz's inequality, we have
\begin{align*}
\sum_{x\in X}\left(\sum_{y\in X}\phi(x,y)\xi(y)\right)^2 &= \sum_{x\in X}\left(\sum_{y\in X}\phi(x,y)^{1/2}\phi(x,y)^{1/2}\xi(y)\right)^2\\
&\leq \sum_{x\in X}\left(\sum_{y\in X}\phi(x,y)\right)\left(\sum_{y\in X}\phi(x,y)\xi(y)^2\right)\\
&\leq\left(\sup_{x\in X}\sum_{y\in X}\phi(x,y)\right)^2\|\xi\|^2_{l^2(X)}~~~~~~~~~~~~~~~~~~~~~~~~~~~~~~~~~~~~~~~~~~~~~~~~~~~~~~~~~~~(\ast)
\end{align*}
On the other hand, by Lemma \ref{croissance}, one can find $B,L>0$ such that $|S(x,n)|\leq BL^n$ for every $n$ and every $x$. Then, if we fix $N\geq n_0$ satisfying $Le^{-N^{\alpha}/k}<1$, for any $x\in X$ we obtain that
\begin{align*}
\sum_{y\in X}\phi(x,y)=\sum_{n\geq 0}\sum_{y\in S(x,n)}\phi(x,y)=\sum_{0\leq n\leq N}\sum_{y\in S(x,n)}\phi(x,y) + \sum_{n> N}\sum_{y\in S(x,n)}\phi(x,y)
\end{align*}
As $\phi(x,y)\leq 1$, the first term is at most $B\sum_{0\leq n\leq N}L^n$. Concerning the second term, by our choice of $N$, we have $\phi(x,y)\leq e^{-\frac{n^{1+\alpha}}{k}}$ for every $n> N$ and $y\in S(x,n)$. Therefore,
$$
\sum_{n> N}\sum_{y\in S(x,n)}\phi(x,y)\leq B\sum_{n> N}L^n e^{-\frac{n^{1+\alpha}}{k}}\leq B\sum_{n> N}\left(L e^{-N^{\alpha}/k}\right)^n ~~~~~~~~~~~~~~~~~~~~~~~~~~~~~~~~~~~~~~~~~~~~~~~~~~~~~~~~~(\ast\ast)
$$
which is the remainder of a convergent geometric series and does not depend on $x$.\\

\noindent Hence $\mathcal{U}_{\phi}:l^2(X)\rightarrow l^2(X)$ is a bounded operator with $1\leq \|\mathcal{U}_{\phi}\|$, as $\mathcal{U}_{\phi}(\delta_x)(x)=1$ for every $x\in X$. Moreover, $\mathcal{U}_{\phi}$ is self-adjoint and positive as $\phi$ is symmetric and positive definite. In particular, one can consider the positive square root of $\mathcal{U}_{\phi}$ which can be represented as (see for instance Theorem VI.9 in \cite{rs})
$$
\mathcal{V}_{\phi}:=\|\mathcal{U}_{\phi}\|^{1/2}\sum_{n\geq 0}a_n\left(I - \frac{\mathcal{U}_{\phi}}{\|\mathcal{U}_{\phi}\|}\right)^n
$$
where $I:=\textrm{Id}_{l^2(X)}$ and $\sqrt{1-z}=\sum_{n\geq 0}a_nz^n$ converges absolutely for every $z$ such that $|z|\leq 1$.\\

\noindent Now fix $M_0$ large enough such that
$$
\left\|\mathcal{V}_{\phi} - \|\mathcal{U}_{\phi}\|^{1/2}\sum_{0\leq n\leq M_0}a_n\left(I - \frac{\mathcal{U}_{\phi}}{\|\mathcal{U}_{\phi}\|}\right)^n\right\| < \frac{\epsilon}{4(4\|\mathcal{U}_{\phi}\|^{1/2} + 1)}.
$$
We fix also $M\geq R$ such that
$$
\sum_{n> M}\left(L e^{-N^{\alpha}k}\right)^n < \frac{\|\mathcal{U}_{\phi}\|^{1/2} \epsilon}{2^{M_0 + 3}(4\|\mathcal{U}_{\phi}\|^{1/2} + 1)B}
$$
Then we denote
$$
\phi_M:X\times X\rightarrow\mathbb{R}~,~ (x,y)\mapsto
\begin{cases}
\phi(x,y)~~~~\textrm{if $d(x,y)\leq M$}\\
0~~~~~~~~~~~~\textrm{otherwise}
\end{cases}
$$
and we define $\mathcal{U}_{\phi_M}:l^2(X)\rightarrow l^2(X)$ by $\mathcal{U}_{\phi_M}(\xi)(x):=\sum_{y\in X}\phi_M(x,y)\xi(y)$ for every $\xi\in l^2(X)$.\\

\noindent If $\phi_M$ was positive definite the proof would be finished but a priori there is no reason for that.\\

\noindent By $(\ast)$, $(\ast\ast)$ and our choice of $M$, we have
$$
\|\mathcal{U}_{\phi} - \mathcal{U}_{\phi_M}\|\leq \frac{\|\mathcal{U}_{\phi}\|^{1/2} \epsilon}{2^{M_0 +3}(4\|\mathcal{U}_{\phi}\|^{1/2} + 1)}~~~~~~~~~~~~~~~~~~~~~~~~~~~~~~~~~~~~~~~~~~~~~~~~~~~~~~~~~~~~~~~~~~~~~~~~~~~~~~~~~ (\ast\ast\ast)
$$

\noindent Finally, by setting
$$
\mathcal{W}:=\|\mathcal{U}_{\phi}\|^{1/2}\sum_{0\leq n \leq M_0}a_n\left(I - \frac{\mathcal{U}_{\phi_M}}{\|\mathcal{U}_{\phi}\|}\right)^n
$$
we can consider the positive definite kernel
$$
\psi:X\times X\rightarrow\mathbb{R}~,~ (x,y)\mapsto \left\langle\mathcal{W}(\delta_x),\mathcal{W}(\delta_y)\right\rangle.
$$

\noindent To prove that $|1-\psi(x,y)| < \epsilon$ whenever $d(x,y)\leq R$, it suffices to show that
 $|\phi(x,y)-\psi(x,y)| < \epsilon/2$ for every $x,y\in X$. We have
\begin{align*}
|\phi(x,y)-\psi(x,y)|&=|\left\langle\mathcal{U}_{\phi}(\delta_x),\delta_y\right\rangle - \left\langle\mathcal{W}^{\ast}\mathcal{W}(\delta_x),\delta_y\right\rangle|\\
&= |\left\langle\left(\mathcal{V}_{\phi}^{\ast}\mathcal{V}_{\phi} - \mathcal{W}^{\ast}\mathcal{W}\right)(\delta_x),\delta_y\right\rangle|\\
&\leq\left\|\mathcal{V}_{\phi}^{\ast}\mathcal{V}_{\phi} - \mathcal{W}^{\ast}\mathcal{W}\right\|\\
&\leq\left(\|\mathcal{V}_{\phi}\| + \|\mathcal{W}\|\right)\left\|\mathcal{V}_{\phi} - \mathcal{W}\right\|\\
&\leq 2\left\|\mathcal{V}_{\phi} - \mathcal{W}\right\|\|\mathcal{V}_{\phi}\| + \left\|\mathcal{V}_{\phi} - \mathcal{W}\right\|^2\\
&\leq \left(2\|\mathcal{V}_{\phi}\| + 1 \right)\left\|\mathcal{V}_{\phi} - \mathcal{W}\right\|\\
& \leq \left(4\|\mathcal{U}_{\phi}\|^{1/2} + 1 \right)\left\|\mathcal{V}_{\phi} - \mathcal{W}\right\|
\end{align*}
But we have
\begin{align*}
\left\|\mathcal{V}_{\phi} - \mathcal{W}\right\| &\leq \left\|\mathcal{V}_{\phi} - \|\mathcal{U}_{\phi}\|^{1/2}\sum_{0\leq n\leq M_0}a_n\left(I - \frac{\mathcal{U}_{\phi}}{\|\mathcal{U}_{\phi}\|}\right)^n\right\| + \left\|\|\mathcal{U}_{\phi}\|^{1/2}\sum_{0\leq n\leq M_0}a_n\left(I - \frac{\mathcal{U}_{\phi}}{\|\mathcal{U}_{\phi}\|}\right)^n - \mathcal{W}\right\|\\
&< \frac{\epsilon}{4(4\|\mathcal{U}_{\phi}\|^{1/2} + 1)} + \left\|\|\mathcal{U}_{\phi}\|^{1/2}\sum_{0\leq n\leq M_0}a_n\left(I - \frac{\mathcal{U}_{\phi}}{\|\mathcal{U}_{\phi}\|}\right)^n - \mathcal{W}\right\|
\end{align*}
and moreover, as $|a_n|\leq 1$ for every $n$, using the inequality $\|A^n - B^n\|\leq 2^n\|A - B\|$ (when $\|A\|\leq 1$ and $\|B\|\leq 2$), we obtain
\begin{align*}
\left\|\|\mathcal{U}_{\phi}\|^{1/2}\sum_{0\leq n\leq M_0}a_n\left(I - \frac{\mathcal{U}_{\phi}}{\|\mathcal{U}_{\phi}\|}\right)^n - \mathcal{W}\right\| & \leq |\mathcal{U}_{\phi}\|^{1/2}\sum_{1\leq n\leq M_0}|a_n|\left\|\left(I - \frac{\mathcal{U}_{\phi}}{\|\mathcal{U}_{\phi}\|}\right)^n - \left(I - \frac{\mathcal{U}_{\phi_M}}{\|\mathcal{U}_{\phi}\|}\right)^n\right\|\\
&\leq 2^{M_0 + 1}\frac{\|\mathcal{U}_{\phi} - \mathcal{U}_{\phi_M}\|}{\|\mathcal{U}_{\phi}\|^{1/2}}\\
& < \frac{\epsilon}{4(4\|\mathcal{U}_{\phi}\|^{1/2} + 1)}
\end{align*}
the last inequality coming from $(\ast\ast\ast)$. Hence we deduce that
$$
\left\|\mathcal{V}_{\phi} - \mathcal{W}\right\|< \frac{\epsilon}{2(4\|\mathcal{U}_{\phi}\|^{1/2} + 1)}
$$
and therefore $|\phi(x,y)-\psi(x,y)| < \epsilon/2$ for every $x,y\in X$.\\

\noindent It remains to show that $\psi$ has finite width (more precisely, it has width at most $2M_0M$). Indeed, if $x,y\in X$ are such that
$$
\psi(x,y)=\sum_{z\in X}\mathcal{W}(\delta_x)(z)\mathcal{W}(\delta_y)(z)\neq 0,
$$
then there exists at least one $z\in X$ such that $\mathcal{W}(\delta_x)(z)\neq 0$ and $\mathcal{W}(\delta_y)(z)\neq 0$. But for every $t\in X$,
$$
\mathcal{W}(\delta_t)(z)=\|\mathcal{U}_{\phi}\|^{1/2}\sum_{s\in X}\sum_{0\leq n \leq M_0}a_n\left(\delta - \frac{\phi_M}{\|\mathcal{U}_{\phi}\|}\right)^{\ast n}(z,s)~\delta_t(s)=\|\mathcal{U}_{\phi}\|^{1/2}\sum_{0\leq n \leq M_0}a_n\left(\delta - \frac{\phi_M}{\|\mathcal{U}_{\phi}\|}\right)^{\ast n}(z,t),
$$
where $\delta(z,z)=1$, $\delta(z,t)=0$ if $z\neq t$, and where $\lambda\ast\mu$ denotes the convolution product, i.e. $\lambda\ast\mu(z,t)=\sum_{s\in X}\lambda(z,s)\mu(s,t)$ ($\nu^{\ast n}$ being the $n$-fold convolution product of $\nu$ with itself). Hence, with these notations, for some $p,q\leq M_0$ we have
$$
\left(\delta - \frac{\phi_M}{\|\mathcal{U}_{\phi}\|}\right)^{\ast p}(z,x)\neq 0~~~~\textrm{and}~~~~\left(\delta - \frac{\phi_M}{\|\mathcal{U}_{\phi}\|}\right)^{\ast q}(z,y)\neq 0
$$
Now, by definition, $\textrm{supp}\left(\delta - \frac{\phi_M}{\|\mathcal{U}_{\phi}\|}\right)\subset\{(a,b)\in X~\mid~ d(a,b)\leq M\}$, and it is easy to show by induction on $l$ that $\textrm{supp}\left(\left(\delta - \frac{\phi_M}{\|\mathcal{U}_{\phi}\|}\right)^{\ast l}\right)\subset\{(a,b)\in X~\mid~ d(a,b)\leq lM\}$. We deduce that $d(x,y)\leq (p+q)M\leq 2M_0M$.

\end{proof}

\section{Groups acting on metric spaces}
\noindent If $X$ is a discrete metric space, by the $\check{\textrm{S}}$varc-Milnor Lemma, any countable group $G$ acting properly and co-compactly by isometries on $X$ is quasi-isometric (thus coarsely equivalent) to $X$. Therefore, property A of $G$ is equivalent to property A of $X$. In this setting, it is easy to make explicit Ozawa kernels:

\begin{prop}\label{noyaux pour actions propres}
Let $G$ be a countable group acting properly by isometries on a discrete metric space $X$. Suppose that there exists $x_0\in X$ such that the orbit $G\cdot x_0$ has property A. Let $R,\epsilon >0$, let $K(R):=\max_{g\in B_G(g,R)}d_X(x_0,gx_0)$ and let $\phi$ be a $(K(R),\epsilon)$-Ozawa kernel on $G\cdot x_0$.  Then
$$
\widetilde{\phi}:G\times G\rightarrow\mathbb{R}~,~ (g,g')\mapsto \phi(gx_0,g'x_0)
$$
defines an $(R,\epsilon)$-Ozawa kernel on $G$ .
\end{prop}

\begin{proof}[Proof of Proposition $\ref{noyaux pour actions propres}$]
By definition, $\phi$ is a positive definite kernel on $G$ satisfying $|1-\widetilde{\phi}(g,g')| < \epsilon $, whenever $d_G(g,g')\leq R$. Moreover, there exists $S >0$ such that $\textrm{supp}(\phi)\subset\{(gx_0,g'x_0)\in G\cdot x_0\times G\cdot x_0 \mid d(gx_0,g'x_0)\leq S\}$.
Hence, for every $g,g'\in G$ such that $\widetilde{\phi}(g,g')\neq 0$, we have $g^{-1}g'x_0\in B_X(x_0,S)\cap G\cdot x_0$ which is finite (by properness). Then, if $B_X(x_0,S)\cap G\cdot x_0=\{g_1\cdot x_0,\ldots,g_N\cdot x_0\}$, we obtain that $\textrm{supp}(\widetilde{\phi})\subset \bigcup_{i=1}^Ng_iG_0$, where
$G_0$ is the (finite) stabilizer of $x_0$. This last set is finite, and if $\widetilde{S}$ denotes the diameter of this set, we have $\textrm{supp}(\widetilde{\phi})\subset\{(g,g')\in G\times G \mid d(g,g')\leq \widetilde{S}\}$.

\end{proof}

\noindent It is well-known that the 0-skeleton (endowed with the induced metric) of a finite dimensional CAT(0) cube complex with bounded geometry has property A (see \cite{cn}), then Proposition \ref{noyaux pour actions propres} gives a partial generalization of Theorem B in \cite{cn}:

\begin{cor}
Any countable group acting properly by isometries on a finite dimensional CAT(0) cube complex with bounded geometry has property A.
\end{cor}

\begin{rem}
Note that according to a recent result of Brodzki, Campbell, Guentner,
Niblo and Wright, the previous corollary remains true without the
assumption of bounded geometry.
\end{rem}

\noindent It was pointed out to me by E. Guentner that a stronger result actually holds. The following theorem is a direct consequence of ideas developed in \cite{dg1} but does not appear explicitly. Here we give a self-contained proof, expliciting Ozawa kernels, in the setting of countable groups.

\begin{thm}\label{noyaux pour actions avec stab ayant A}
Let $G$ be a countable group acting by isometries on a discrete metric space with bounded geometry $X$. Assume that there exists $x_0\in X$ such that both the orbit $G\cdot x_0$ and the stabilizer $G_0$ of $x_0$ have property A. Then $G$ has property A.
\end{thm}

\begin{rem}
In contrast to what happens in the case of a proper action, in general the kernel $\widetilde{\phi}$ in the proof of Proposition $\ref{noyaux pour actions propres}$ will not satisfy the support condition. But, when we have property A (instead of finiteness) for $G_0$, one can nevertheless use the following idea: if we cover a metric space (with bounded geometry) with subsets having property A and if we have an adapted subordinated partition of unity to glue them together, one obtains an Ozawa kernel on the whole metric space.

\end{rem}

\begin{proof}[Proof of Theorem $\ref{noyaux pour actions avec stab ayant A}$]
Let $R,\epsilon>0$. We begin by proceeding as in the proof of Proposition $\ref{noyaux pour actions propres}$. Denote by $K(R):=\max_{g\in B_G(g,R)}d_X(x_0,gx_0)$, and let $\phi_0$ be a $(K(R),\frac{\epsilon^2}{8})$-Ozawa kernel for $G\cdot x_0$ given by Proposition $\ref{forme des noyaux}$. That is,
$\phi_0:(gx_0,g'x_0)\mapsto \langle\lambda(gx_0),\lambda(g'x_0)\rangle$ for some $\lambda:G\cdot x_0\rightarrow l^2(G\cdot x_0)_1$ with $\textrm{supp}(\lambda(gx_0))\subset B_X(gx_0,S_0)$ for every $x$ and $\lambda(gx_0)(g'x_0)\geq 0$ for every $g$. For every $g\in G$, we define
$$
U_g:=\{g'\in G ~\mid~ gx_0\in\textrm{supp}(\lambda(g'x_0))\}
$$
and $\alpha_g:G\rightarrow\mathbb{R},g'\mapsto\lambda(g'x_0)^2(gx_0)$. Hence we obtain a cover $\mathcal{U}_G:=\{U_g\}_{g\in G}$ of $G$ and a partition of unity $\alpha:=\{\alpha_g\}_{g\in G}$ subordinated to it (i.e., $\alpha_g(g')=0$ if $g'\notin U_g$ and $\sum_{g\in G}\alpha_g(g')=1$ for every $g'\in G$). Moreover, if $d_G(g_1,g_2)\leq R$, by Cauchy-Schwarz's inequality,
\begin{align*}
\sum_{g\in G}|\alpha_g(g_1)-\alpha_g(g_2)|&=\left\|\lambda(g_1x_0)^2-\lambda(g_2x_0)^2\right\|_{l^1(G\cdot x_0)}\\
&\leq 2 \left\|\lambda(g_1x_0)-\lambda(g_2x_0)\right\|_{l^2(G\cdot x_0)}\\
&=2\sqrt{2-2\phi_0(g_1x_0,g_2x_0)}\leq\epsilon.
\end{align*}
For every $g\in G$, if $B_X(x_0,S_0)\cap G\cdot x_0=\{g_1\cdot x_0,\ldots,g_N\cdot x_0\}$ (this set is finite by the bounded geometry condition on $X$), we obtain that $U_g\subset \bigcup_{i=1}^Ngg_iG_0$. For convenience, we denote
$$
X_g^i:=gg_iG_0~,~~~ X_g:=\bigcup_{i=1}^NX_g^i.
$$
For any subset $A$ of $X_g$ and for any $L>0$, we denote by $A(L):=\{x\in X_g ~\mid~d_G(x,A)\leq L\}$ the $L$-neighborhood of $A$ in $X_g$.\\

\noindent Now the main part of the proof is to construct an Ozawa kernel on each $X_g$ (then on each $U_g$) using Ozawa's property on $G_0$.\\

\noindent Let $R_1:=3R$, $\epsilon_1:=\frac{\epsilon^2}{4}$, and let $L\geq\frac{4N(2N+1)R_1}{\epsilon_1}$. Then $X_g^i(L)$ is a finite cover of $X_g$ of multiplicity $\leq N$ and with Lebesgue number $\geq L$ (see Definition $\ref{asdim}$). Then defining for $x\in X_g$,
$$
\delta_g^i(x):=\frac{d_G\left(x,X_g\smallsetminus X_g^i(L)\right)}{\sum_{j=1}^Nd_G\left(x,X_g\smallsetminus X_g^j(L)\right)},
$$
we obtain a partition of unity subordinated to the cover $\{X_g^i(L)\}_{i=1}^N$ satisfying, by Lemma $\ref{noyau de recouvrement}$ and by our choice for $L$,
$$
\sum_{i=1}^N|\delta_g^i(x)-\delta_g^i(y)| \leq \frac{\epsilon_1}{2},
$$
whenever $d_G(x,y)\leq R_1$.\\

\noindent Now, we fix a $(2(L+R_1),\frac{\epsilon_1^2}{32})$-Ozawa kernel on $G_0$, $\psi_0:(z,z')\mapsto \langle\mu(z),\mu(z')\rangle$, for some $\mu:G_0\rightarrow l^2(G_0)_1$ with $\textrm{supp}(\mu(z))\subset B_G(z,S_1)$ for every $z$. Hence we deduce a $(2(L+R_1),\frac{\epsilon_1^2}{32})$-Ozawa kernel (with the same support) on each $X_g^i$ by setting
$$
\psi_i(gg_iz,gg_iz'):=\psi_0(z,z')
$$
for every $z,z'\in G_0$. In order to obtain an appropriate Ozawa kernel on $X_g$, for any $z\in X_g$, we fix $p_i(z)\in X_g^i$ such that
$d_G(z,X_g^i)=d_G(z,p_i(z))$. Then we put $\widetilde{X_g}:=\coprod_{i=1}^N\{i\}\times X_g^i$, and we define
$\sigma_g:X_g\rightarrow l^2(\widetilde{X_g})~,~z\mapsto\sigma_g(z)$,
where
$$
\sigma_g(z)(i,x):=
\delta_g^i(z)^{1/2}\mu((gg_i)^{-1}p_i(z))((gg_i)^{-1}x).
$$
Observe that if $\sigma_g(z)(i,x)\neq 0$ for some $x\in X_g^i$, then $\mu((gg_i)^{-1}p_i(z))((gg_i)^{-1}x)\neq 0$, $z\in X_g^i(L)$ and we have $d_G((gg_i)^{-1}p_i(z),(gg_i)^{-1}x)=d_G(x,p_i(z))\leq S_1$. Therefore,
$$
d_G(x,z)\leq d_G(x,p_i(z)) + d_G(p_i(z),z)\leq S_1+L~~~~~~~~~~~~~~~~~~~~~~~~~~~~~~~~~~~~~~~~~~~~~~~~~~~~~~~~~~~~~~~~~~~~~~~~(\ast)
$$
Moreover, for $z_1,z_2\in X_g$ such that $d_G(z_1,z_2)\leq R_1$, one has

\begin{align*}
\|\sigma_g(z_1)-\sigma_g(z_2)\|_{l^2(\widetilde{X_g})}^2&=\sum_{i=1}^N\sum_{x\in G_0}\left|\delta_g^i(z_1)^{1/2}\mu((gg_i)^{-1}p_i(z_1))(x) - \delta_g^i(z_2)^{1/2}\mu((gg_i)^{-1}p_i(z_2))(x)\right|^2\\
&\leq \sum_{i=1}^N\sum_{x\in G_0}\left|\delta_g^i(z_1)\mu((gg_i)^{-1}p_i(z_1))^2(x) - \delta_g^i(z_2)\mu((gg_i)^{-1}p_i(z_2))^2(x)\right|\\
&\leq \underbrace{\left(\sum_{i=1}^N\delta_g^i(z_1)\right)\left(\sum_{x\in G_0}\left|\mu((gg_i)^{-1}p_i(z_1))^2(x) - \mu((gg_i)^{-1}p_i(z_2))^2(x)\right|\right)}_{(I)}\\
&+ \underbrace{\left(\sum_{x\in G_0}\mu((gg_i)^{-1}p_i(z_2))^2(x)\right)\left(\sum_{i=1}^N\left|\delta_g^i(z_1) - \delta_g^i(z_2)\right|\right)}_{(II)}.
\end{align*}
However, on one hand, $\sum_{i=1}^N\delta_g^i(z_1)=1$ and thus there exists $i$ such that $z_1\in X_g^i(L)$. In this case, we have $z_2\in (X_g^i(L))(R_1)$, then $d_G(p_i(z_1),p_i(z_2))\leq 2(L+R_1)$. By Cauchy-Schwarz's inequality,
\begin{align*}
\sum_{x\in G_0}\left|\mu((gg_i)^{-1}p_i(z_1))^2(x) - \mu((gg_i)^{-1}p_i(z_2))^2(x)\right|&=\left\|\mu((gg_i)^{-1}p_i(z_1))^2 - \mu((gg_i)^{-1}p_i(z_2))^2\right\|_{l^1(G_0)}\\
&\leq 2\left\|\mu((gg_i)^{-1}p_i(z_1)) - \mu((gg_i)^{-1}p_i(z_2))\right\|_{l^2(G_0)}\\
&=2\sqrt{2}\sqrt{1-\psi_0((gg_i)^{-1}p_i(z_1),(gg_i)^{-1}p_i(z_2))}\\
&\leq \frac{\epsilon_1}{2}.
\end{align*}
We conclude that $(I)\leq \frac{\epsilon_1}{2}$. On the other hand,
$$
\sum_{x\in G_0}\mu((gg_i)^{-1}p_i(z_2))^2(x)=\left\|\mu((gg_i)^{-1}p_i(z_2))\right\|^2_{l^2(G_0)}=1.
$$
Hence, part $(II)$ is less or equal to
$$
\sum_{i=1}^N\left|\delta_g^i(z_1) - \delta_g^i(z_2)\right|\leq \frac{\epsilon_1}{2}.
$$
It follows that $\left\|\sigma_g(z_1)-\sigma_g(z_2)\right\|_{l^2(\widetilde{X_g})}^2\leq \epsilon_1$ whenever $d_G(z_1,z_2)\leq R_1$.\\

\noindent Now let us deduce an Ozawa kernel on $U_g$. For every $x\in X_g$ let us fix $q_g(x)\in U_g$ such that
$d_G(x,q_g(x))=d_G(x,U_g)$ and define $\tau_g:U_g\rightarrow l^2(U_g)~,~u\mapsto\tau_g(u)$, where
$$
\tau_g(u):w\mapsto\left(\sum_{\substack{(i,x)\in \widetilde{X_g}:\\q_g(x)=w}}\sigma_g(u)^2(i,x)\right)^{1/2}.
$$
In particular, as $\bigcup_{w\in U_g}\{(i,x)\in \widetilde{X_g}~\mid~q_g(x)=w\}=\widetilde{X_g}$, we have $\left\|\tau_g(u)\right\|_{l^2(U_g)}=\left\|\sigma_g(u)\right\|_{l^2(\widetilde{X_g})}=1$. By Minkowski's inequality, $d_G(u,u')\leq R_1$ implies
$$
\left\|\tau_g(u)-\tau_g(u')\right\|_{l^2(U_g)}^2\leq\left\|\sigma_g(u)-\sigma_g(u')\right\|_{l^2(\widetilde{X_g})}^2\leq\epsilon_1~~~~~~~~~~~~~~~~~~~~~~~~~~~~~~~~~~~~~~~~~~~~~~~~~~~~~~~~~~~~~~~~~~(\ast\ast)
$$
Moreover, if $\tau_g(u)(w)\neq 0$, there exists $(i,x)\in \widetilde{X_g}$ such that $q_g(x)=w$ with $\sigma_g(u)(i,x)\neq 0$. Thus, by $(\ast)$,
$d_G(u,x)\leq S_1+L$ and
$$
d_G(u,w)\leq d_G(u,x) + d_G(x,q_g(x))\leq 2d_G(u,x)\leq 2(S_1+L).
$$
Then $\textrm{supp}(\tau_g(u))\subset B_G(u,2(S_1+L))$ for every $u$.\\

\noindent Now, for every $g\in G$ and for every $x\in U_g(R)$, we fix $r_g(x)\in U_g$ such that $d_G(x,r_g(x))=d_G(x,U_g)\leq R$. Then we set
$$
\nu_g:U_g(R)\rightarrow l^2(U_g)~,~ x\mapsto\nu_g(x):=\tau_g(r_g(x)).
$$
By definition, if $d_G(x,z)\leq R$ we have $d_G(r_g(x),r_g(z))\leq 3R=R_1$, hence by $(\ast\ast)$
$$
\left\|\nu_g(x) - \nu_g(z)\right\|_{l^2(U_g)}^2=\left\|\tau_g(r_g(x)) - \tau_g(r_g(z))\right\|_{l^2(U_g)}^2\leq \epsilon_1
$$
and if $\nu_g(x)(u):=\tau_g(r_g(x))(u)\neq 0$, we have $d_G(u,r_g(x))\leq 2(S_1+L)$ and $d_G(u,x)\leq R + 2(S_1+L)$, i.e., $\textrm{supp}(\nu_g(x))\subset B_G(x,2(S_1+L)+R)$ for every $x$.\\

\noindent Finally, we extend each $\nu_g$ on $G$ by 0 outside $U_g(R)$, we put $\mathcal{U}_G:=\coprod_{g\in G}\{g\}\times U_g$ and we consider the positive definite kernel
$$
\psi:G\times G\rightarrow\mathbb{R}~,~(g_1,g_2)\mapsto\langle\kappa(g_1),\kappa(g_2)\rangle=\sum_{g\in G}\sum_{x\in U_g}\kappa(g_1)(g,x)\kappa(g_2)(g,x)
$$
where $\kappa:G\rightarrow l^2(\mathcal{U}_G)$ is defined by
$$
\kappa(g')(g,x):=\alpha_g^{1/2}(g')\nu_g(g')(x).
$$
On the one hand, if $\psi(g_1,g_2)\neq 0$ there exists $(g,x)\in\mathcal{U}_G$ such that $\kappa(g_1)(g,x):=\alpha_g^{1/2}(g_1)\nu_g(g_1)(x)\neq 0$ and $\kappa(g_2)(g,x):=\alpha_g^{1/2}(g_2)\nu_g(g_2)(x)\neq 0$ , i.e. $\nu_g(g_1)(x)\neq 0$, $\nu_g(g_2)(x)\neq 0$ and $g_1,g_2\in U_g$. Then $d_G(x,g_1)\leq 2(S_1+L)+R$ and $d_G(x,g_2)\leq 2(S_1+L)+R$. Hence
$d_G(g_1,g_2)\leq 4(S_1+L)+2R$. Therefore,
$$
\textrm{supp}(\psi)\subset\{(g_1,g_2)\in G\times G~\mid~d_G(g_1,g_2)\leq 4(S_1+L)+2R\}.
$$
On the other hand,
\begin{align*}
\left\|\kappa(g_1)-\kappa(g_2)\right\|_{l^2(\mathcal{U}_G)}^2&=\sum_{g\in G}\sum_{x\in U_g}\left|\kappa(g_1)(g,x)-\kappa(g_2)(g,x)\right|^2\\
&=\sum_{g\in G}\sum_{x\in U_g}\left|\alpha_g^{1/2}(g_1)\nu_g(g_1)(x)-\alpha_g^{1/2}(g_2)\nu_g(g_2)(x)\right|^2\\
&\leq \sum_{g\in G}\sum_{x\in U_g}\left|\alpha_g(g_1)\nu_g(g_1)^2(x)-\alpha_g(g_2)\nu_g(g_2)^2(x)\right|\\
&\leq \underbrace{\left(\sum_{g\in G}\alpha_g(g_1)\right)\left(\sum_{x\in U_g}\left|\nu_g(g_1)^2(x)-\nu_g(g_2)^2(x)\right|\right)}_{(A)}\\
& + \underbrace{\left(\sum_{x\in U_g}\nu_g(g_2)^2(x)\right)\left(\sum_{g\in G}\left|\alpha_g(g_1)-\alpha_g(g_2)\right|\right)}_{(B)}.
\end{align*}
If $d_G(g_1,g_2)\leq R$, $\sum_{g\in G}\alpha_g(g_1)=1$ and there exists at least one $g\in G$ such that $g_1\in U_g$, then $g_2\in U_g(R)$. Thus, by Cauchy-Schwarz's inequality
\begin{align*}
\sum_{x\in U_g}\left|\nu_g(g_1)^2(x)-\nu_g(g_2)^2(x)\right|&=\left\|\nu_g(g_1)^2-\nu_g(g_2)^2\right\|_{l^1(U_g)}\\
&\leq 2\left\|\nu_g(g_1)-\nu_g(g_2)\right\|_{l^2(U_g)}\\
&\leq 2\sqrt{\epsilon_1}=\epsilon
\end{align*}
i.e. $(A)\leq \epsilon$. For the term $(B)$, one has $\sum_{x\in U_g}\nu_g(g_2)^2(x)=0$ if $g_2\notin U_g(R)$ and $\sum_{x\in U_g}\nu_g(g_2)^2(x)=\|\nu_g^2(g_2)\|^2_{l^2(U_g)}=1$ if $g_2\in U_g(R)$. Therefore in all cases
$$
(B)\leq \sum_{g\in G}\left|\alpha_g(g_1)-\alpha_g(g_2)\right|\leq \epsilon
$$
hence $1 - \psi(g_1,g_2)=\frac{1}{2}\left\|\kappa(g_1)-\kappa(g_2)\right\|_{l^2(\mathcal{U}_G)}^2\leq\epsilon$. This concludes the proof.\\

\end{proof}

\bigskip

\noindent What precedes can be summarized as follows:

\begin{prop}\label{noyaux pour actions}
Let $G$ be a countable group acting by isometries on a discrete metric space with bounded geometry $X$. Assume that there exists $x_0\in X$ such that both the orbit $G\cdot x_0$ and the stabilizer $G_0$ of $x_0$ have property A. Let $R,\epsilon>0$. Let $K(R):=\max_{g\in B_G(g,R)}d_X(x_0,gx_0)$ and fix an $(K(R),\frac{\epsilon^2}{8})$-Ozawa kernel for $G\cdot x_0$,
$\phi_0:(gx_0,g'x_0)\mapsto \langle\lambda(gx_0),\lambda(g'x_0)\rangle$, for some $\lambda:G\cdot x_0\rightarrow l^2(G\cdot x_0)_1$ with $\textrm{supp}(\lambda(gx_0))\subset B_X(gx_0,S_0)$ for every $x$ (and $\lambda(gx_0)(g'x_0)\geq 0$ for every $g$). For every $g\in G$ we define
$$
U_g:=\{g'\in G ~\mid~ gx_0\in\textrm{supp}(\lambda(g'x_0))\}
$$
If $B_X(x_0,S_0)\cap G\cdot x_0=\{g_1\cdot x_0,\ldots,g_N\cdot x_0\}$, for each $g\in G$ we have $U_g\subset \bigcup_{i=1}^Ngg_iG_0:=X_g$.
\noindent Let $L\geq\frac{48N(2N+1)R}{\epsilon^2}$, fix a $(2(L+3R),\frac{\epsilon^4}{512})$-Ozawa kernel on $G_0$, $\psi_0:(z,z')\mapsto \langle\mu(z),\mu(z')\rangle$, for some $\mu:G_0\rightarrow l^2(G_0)_1$ with $\textrm{supp}(\mu(z))\subset B_G(z,S_1)$ for every $z$. Finally set for every $x\in X_g$
$$
\delta_g^i(x):=\frac{d_G\left(x,X_g\smallsetminus (gg_iG_0)(L)\right)}{\sum_{j=1}^Nd_G\left(x,X_g\smallsetminus (gg_iG_0)(L)\right)}
$$
Hence
\begin{align*}
&\zeta_g(u,u'):=\\
&\sum_{w\in U_g}\left[\sum_{\substack{(i,x)\in \widetilde{X_g}:\\q_g(x)=w}}\delta_g^i(u)~\mu((gg_i)^{-1}p_i(u))^2((gg_i)^{-1}x)\right]^{1/2}
\left[\sum_{\substack{(i,x)\in \widetilde{X_g}:\\q_g(x)=w}}\delta_g^i(u')~\mu(gg_i)^{-1}p_i(u'))^2((gg_i)^{-1}x)\right]^{1/2}
\end{align*}
defines an Ozawa kernel on $U_g$. Here $\widetilde{X_g}:=\coprod_{i=1}^N\{i\}\times gg_iG_0$ and $p_i$, $q_g$ satisfy $p_i(z)\in X_g^i$, $q_g(z)\in U_g$,
$d_G(z,p_i(z))=d_G(z,gg_iG_0)$ and $d_G(z,q_g(z))=d_G(z,U_g)$ for any $z\in X_g$.\\

\noindent Therefore
$$
\psi(g_1,g_2):=\sum_{g\in G}\text{{\Large$\chi$}}_{U_g(R)}(g_1)~\text{{\Large$\chi$}}_{U_g(R)}(g_2)~\lambda(g_1x_0)(gx_0)~\lambda(g_2x_0)(gx_0)~
\zeta_g(r_g(g_1),r_g(g_2))~~~~~~~~~~~~~~~~~~~~~~~~~~~~ (\sharp)
$$
is an $(R,\epsilon)$-Ozawa kernel on $G$ with $\textrm{supp}(\psi)\subset\{(g_1,g_2)\in G\times G~\mid~d_G(g_1,g_2)\leq 4(S_1+L)+2R\}$, where  $r_g(x)\in U_g$ and $d_G(x,r_g(x))=d_G(x,U_g)$ for every $x\in U_g(R)$.

\end{prop}

\bigskip

\begin{cor}\label{groupe ayant sous groupe et quotient avec A}
Let $G$ be a countable group and $H$ a subgroup with property A. If the set of left cosets $G/H$ (endowed with the quotient metric) has property A, then $G$ has property A.
\end{cor}

\noindent Note that we recover (in a rather complicated way) Theorem $\ref{thm extensions}$, and, moreover, we obtain an alternative proof of the stability of property A under taking certain amalgamated free products and HNN-extensions (see \cite{dg} for the general case):

\begin{cor}\label{produits libres et extensions HNN quand les indices sont finis}
Let $G$ and $H$ be two countable groups.
\begin{enumerate}
\item[$(i)$] If $G$ and $H$ have property A, and if $K$ is a common subgroup of finite index both in $G$ and $H$, then $G\ast_K H$ has property A.

\item[$(ii)$] Assume that $G$ has property A, and that $H$ is a finite index subgroup of $G$. Let $\theta:H\rightarrow G$ be a monomorphism such that $\theta(H)$ also have finite index in $G$ , then the HNN-extension $\textrm{HNN}(G,H,\theta)$ has property A.
\end{enumerate}
\end{cor}

\begin{proof}[Proof of Corollary $\ref{produits libres et extensions HNN quand les indices sont finis}$]
The finite index hypothesis ensures bounded geometry of the corresponding Bass-Serre trees (endowed with the simplicial metric) on which the groups act by isometries. The stabilizers of the vertex are all isometric to $G$ (in the case of HNN-extensions), or, to $G$ or $H$ (in the case of amalgamated free products). Hence, it suffices to use property A for trees (which can be shown exactly in the same way as for free groups) and to apply Theorem $\ref{noyaux pour actions avec stab ayant A}$.
\end{proof}

\section{Applications}

\subsection{Hyperbolic groups}

It is a result due to J. Roe (see \cite{roe}) that finitely generated hyperbolic groups have finite asymptotic dimension. Explicit covers given in \cite{roe} allow us to exhibit explicit Ozawa kernels. Let $\Gamma=\left<S\right>$ be a finitely generated $\delta$-hyperbolic group (endowed with the length function $|\cdot|_S$ associated to $S$), and let $(\cdot\mid\cdot)$ denote the Gromov product on $\Gamma$. Let $N_{\delta}$ denote the number such that each ball of radius $R+6\delta$ can be covered by at most $N_{\delta}$ balls of radius $R$. Let $\epsilon,R>0$, and $L \geq \frac{(2N_{\delta}+1)(4N_{\delta}+3)R}{\epsilon}$. For each $k\geq 1$, we fix a maximal subset $\{\gamma_{ik}\}_{i=1}^{n_k}$ in the sphere of radius $kL$ such that $d_{\Gamma}(\gamma_{ik},\gamma_{jk})> L$ for every $i\neq j$. Then we define
$$
U_{ik}:=\{\gamma\in\Gamma~\mid~ kL\leq|\gamma|_S\leq(k+1)L~,~(\gamma\mid\gamma_{ik})\geq (k-\frac{1}{2})L-\delta\}
$$
This is shown in \cite{roe} that for every $k$, $\{U_{ik}\}_{i=1}^{n_k}$ is a uniformly bounded cover of the annulus
$\{\gamma\in\Gamma~\mid~kL\leq|\gamma|_S\leq(k+1)L\}$ such that every ball of radius less than $L$ meets at most $N_{\delta}$ elements of this cover. Hence, defining $\mathcal{U}:=\bigcup_{k\in\mathbb{N}}\{U_{ik}(L)\}_{i=1}^{n_k}$ (where $U_{ik}(L)$ denote the $L$-neighborhood of $U_{ik}$), it is easy to see that $\mathcal{U}$ is a uniformly bounded cover of $\Gamma$ with multiplicity less than $2N_{\delta}$ and such that $L(\mathcal{U})\geq L$. Therefore $\textrm{asdim}\Gamma\leq 2N_{\delta}-1$, hence, by Theorem \ref{noyaux pour groupes avec asdim finie}, the kernel
$$
\psi:\Gamma\times \Gamma\rightarrow\mathbb{R}~,~(\gamma,\gamma')\mapsto\sum_{i,k}\left(\frac{d(\gamma,X\smallsetminus U_{ik})}{\sum_{j,l}d(\gamma,X\smallsetminus U_{jl})}\right)^{1/2}\left(\frac{d(\gamma',X\smallsetminus U_{ik})}{\sum_{j,l}d(\gamma',X\smallsetminus U_{jl})}\right)^{1/2}
$$
 is an $(R,\epsilon)$-Ozawa kernel for $\Gamma$.

\subsection{CAT(0) cubical groups}

Let $X$ denote a CAT(0) cube complex with bounded geometry and let $X^{(0)}$ denote its 0-skeleton endowed with the path metric on the 1-skeleton of $X$ (when $X$ is finite dimensional this metric is equivalent to the metric induced by $X$). Using Theorem $\ref{noyaux quand compression > 1/2}$ and ideas developed in \cite{cn}, we deduce Ozawa kernels on $X^{(0)}$.  Fix a basepoint $v\in X^{(0)}$ and denote by $H$ the set of ``hyperplanes" in $X$. For every $s\in X^{(0)}$ there is a unique ``normal cube path" $\{C_1,\ldots,C_n\}$. Let us define $w_s:H\rightarrow\mathbb{N}$ by $w_s(h)=i+1$ if $h$ intersects the cube $C_i$ and $w_s(h)=0$ otherwise. Then for every $0<\alpha <1/2$, $f_{\alpha}:X^{(0)}\rightarrow l^2(H)~,~s\mapsto\sum_{h\in H}w_s(h)^{\alpha}\delta_h$ is a uniform embedding with $\rho_+(f_{\alpha})$ linear and $\rho_-(f_{\alpha})(r)=r^{1/2+\alpha}$ for $r$ large enough. Hence, the proof of Theorem $\ref{noyaux quand compression > 1/2}$ provides Ozawa kernels on $X^{(0)}$. More precisely, for a fixed $\alpha>0$ and for $M$ large enough, setting
$$
\phi:X^{(0)}\times X^{(0)}\rightarrow\mathbb{R}~,~ (x,y)\mapsto \exp\left(-\frac{\|f_{\alpha}(x) - f_{\alpha}(y)\|_{l^2(H)}^2}{M}\right)
$$
with the notations in the proof of Theorem \ref{noyaux quand compression > 1/2}, we obtain that
$$
\psi:(x,y)\mapsto \|\mathcal{U}_{\phi}\|\sum_{z\in X^{(0)}}\sum_{0\leq m,n \leq M}a_na_m\left(\delta - \frac{\phi_M}{\|\mathcal{U}_{\phi}\|}\right)^{\ast m}(z,x)\left(\delta - \frac{\phi_M}{\|\mathcal{U}_{\phi}\|}\right)^{\ast n}(z,y)
$$
is an Ozawa kernel on $X^{(0)}$. Now, let us fix a vertex $v_0\in X^{(0)}$. By Proposition \ref{noyaux pour actions propres}, for every group $G$ acting properly (and co-compactly) by isometries on $X$, the kernel $\psi_G:G\times G\rightarrow\mathbb{R}~,~ (g,h)\mapsto\psi(gv_0,hv_0)$ defines an Ozawa kernel on $G$.
\subsection{Baumslag-Solitar groups}

Let $p,q\geq 1$, and let $\textrm{BS}(p,q):=\left< a,b\mid ab^pa^{-1}=b^q\right>$ be a standard presentation of Baumslag-Solitar group. It is the HNN-extension $\textrm{HNN}(G,H,\theta)$, where $G=\left< b \right>\backsimeq\mathbb{Z}$, $H=\left< b^q \right>\backsimeq q\mathbb{Z}$ and $\theta:H\rightarrow G~,~b^q\mapsto b^p$. The group $\textrm{BS}(p,q)$ acts transitively by isometries on its Bass-Serre tree $T_{p,q}$ (which is $(p+q)$-regular) and all the stabilizers of the vertex are isometric to $\mathbb{Z}$. Hence Ozawa kernels on $\textrm{BS}(p,q)$ are given by formula $(\sharp)$ in Proposition \ref{noyaux pour actions} with (for $k,l$ large enough)
$$
\lambda:T_{p,q}\rightarrow l^2(T_{p,q})~,~v\mapsto\frac{\text{{\large$\chi$}}_{A_v}}{k+1}~,~~~~~~ \mu:G\rightarrow l^2(G)~,~x\mapsto\frac{\text{{\large$\chi$}}_{xB_G(e,l)}}{2l+1},
$$
where $A_v$ denotes the intersection of $B_{T_{p,q}}(v,k)$ with the unique geodesic ray starting form $v$ and intersecting a fixed geodesic ray in $T_{p,q}$ as a geodesic ray.

\end{document}